\documentclass[12pt,a4,oneside,reqno]{amsart}

\usepackage[applemac]{inputenc}
\usepackage[french]{babel}
\usepackage{latexsym,enumerate} 
\usepackage{amsmath,amsthm,amsopn,amstext,amscd,amsfonts,amssymb,euler,amsbsy}
\usepackage{mathrsfs}
\usepackage{fullpage}
\usepackage{graphicx}

\newtheorem{proposition}{Proposition}
\newtheorem{lemma}{Lemma}
\newtheorem{theorem}{Theorem}

\newtheorem{rem}{Remark}

\let\epsilon=\varepsilon
\let\eps=\epsilon
\let\phi=\varphi

\let\si=\sigma  
  \let\ka=\kappa
\let\tilde=\widetilde

\newcommand{\field}[1]{\mathbb{#1}}
\newcommand{\R}{\field{R}}
\newcommand{\Rd}{\field{R}^d}

\renewcommand{\L}{\field{L}}

\newcommand{\M}{{\mathcal M}}

\newcommand{\beqn}{\begin{equation}}
\newcommand{\eeqn}{\end{equation}}
\newcommand\eref[1]{(\ref{#1})}

\newcommand{\norm}[1]{\|#1\|}

\newcommand{\PE}[1]{{\lfloor#1\rfloor}}

\def\E{{\field E}}
\def\I{{\field I}}
\def\P{{\field P}}
\def\X{{\field X}}
\def\Mnp{{M_{np}}}

\def\cM{{\field M}}
\def\cH{{\field H}}
\def\lan{\lambda_n}
\def\lam{\lambda}

\def\rhon{{\rho^{\otimes n}}}
\def\rhom{{\rho^{\otimes n}}}
\def\bjk{\beta_{jk}}
\def\hbjk{\hat\beta_{jk}}
\def\tbjk{\tilde\beta_{jk}}
\def\sumin{\sum_{i=1}^n}
\def\sumjk{\sum_{j,k}}
\def\psijk{\psi_{jk}}
\def\hf{\hat f}
\def\f{f_{\rho}}
\def\Ijk{I_{jk}}
\def\bz{{\bf z}}

\title{Thresholding in Learning theory}
\author{Gérard Kerkyacharian and Dominique Picard }
\address{CNRS LPMA, 175 rue du Chevaleret, 75013 Paris, France.  Universit\'e Paris X, 200 Avenue de la R\'epublique 92001 Nanterre cedex.
Universit\'e Paris VII,
175 rue du Chevaleret, 75013 Paris, France }
\begin{document}
\maketitle
\begin{abstract}
In this paper we investigate the problem of learning an unknown bounded function. We  be  emphasize special  cases where it is possible to provide very simple (in terms of computation) estimates  enjoying  in addition the property of being universal : their construction does not depend on a priori knowledge on regularity conditions on the unknown object and still they have almost optimal properties for a whole bunch of functions spaces.
 These estimates are constructed using a thresholding schema, which has proven in the last decade in statistics to have very good properties for recovering  signals with inhomogeneous smoothness but has not been extensively developed in Learning Theory.
 
We will basically consider two particular situations. In the first case, we consider the RKHS situation. In this case, we produce a new  algorithm and investigate   its performances in $L_2(\hat\rho_X)$. The exponential rates of convergences are proved to be almost optimal,  and the regularity assumptions are expressed in simple terms. 
 
 The second case considers a more specified situation where the $X_i$'s are one dimensional and the estimator is a wavelet thresholding estimate. The results are comparable in this setting to those obtained in the RKHS situation as concern the critical value and the exponential rates. The advantage here is that we are able to state the results in the $L_2(\rho_X)$ norm and the regularity conditions are expressed in terms of standard H\"older spaces.

\end{abstract}

\section{Introduction}

In this paper, we  are interested in the problem of learning an unknown function defined on a set $\X$ which takes values
in a set $Y$.  We assume that $\X$ is a compact domain in $\Rd$ and $Y=[-M/2,M/2]$ is a finite interval in $\R$.     This problem, also called regression problem, has a long history in Statistics (many references can be found, for example,  in the following books \cite{IH},  \cite{VG} and \cite{4Authors} ). It has recently drawn much attention in the work of   \cite{CS} and amplified upon in
 \cite{PS}. 

We will assume to observe an $n$ sample $Z_1,\ldots,Z_n$ of $Z=(X,Y)$.
The distribution of $Z$ in denoted by $\rho$.
Our aim is to recover the function $f_\rho$:
$$f_\rho(x)=\E_\rho[Y|X=x].$$

  We shall have as our goal to obtain
estimations to  $f_\rho$ with the error measured in the $L_2(\X,\rho_X)$ norm, or $L_2(\X,\hat \rho_X)$ where $\hat\rho$ is the empirical measure calculated on the $X_i$'s.
  $$\norm{g}_\rho^2=\int_{\X}g(x)^2d\rho(x),\;   \norm{g}_{\hat\rho_X}^2=\frac1n\sumin g(X_i)^2$$
  Given any $\eta>0$, if $\hat f$ is an estimator of $\f$ (i.e. a measurable function of $Z_1,\ldots,Z_n$, taking its values in the set, say, of bounded functions),
\beqn
\label{acc}
\rhom\{\bz:\|\hf-f_\rho\|>\eta\}
\eeqn
measures the confidence we have that the estimator $\hf $ is accurate to tolerance $\eta$.

Contrary to Statistics, where people are mainly concerned with evaluation of moments of 
$\|\hf-f_\rho\|$ (except rare examples, see \cite{K} or \cite{KS}...) Learning Theory focuses on investigating  
 the decay of \eref{acc} as $n\to \infty$ and $\eta$ increases.

Another difference with the Statistics point of view is that one mail goal in Learning Theory is to obtain results with almost no assumptions on the distribution $\rho$. However, it is known that it is not possible to have fast rates of convergence without assumptions and  a large portion of Statistics and Learning Theory proceeds under the
condition that $f_\rho$ is in a known set $\Theta$.    Typical choices of $\Theta$ are compact sets determined by some smoothness condition or by some prescribed rate of decay for a specific approximation process.      Given our prior $\Theta $ and the associated class $\cM(\Theta)$ of measures $\rho$,
it has been defined in \cite{dkpt}, for each $\eta>0$ the {\it accuracy confidence function}
\beqn
\label{accf}
{\bf AC}_n(\Theta,\hf,\eta):=\sup_{\rho\in\cM(\Theta)}\rhom\{\bz:\|f_\rho-\hf\|>\eta\}.
\eeqn
  This quantity  measures a uniform confidence (over the space $\cM(\Theta)$) we have that the estimator $f_\bz$ is accurate to tolerance $\eta$. 
\\
Upper and lower bounds for ${\bf AC}$  have been proved in \cite{dkpt}.
%
%\beqn
%\label{acc2}
%{\bf AC}_m(\Theta,\delta)\ge C\sqrt {\bar N(\Theta,\delta)}e^{-cm\delta^2},
%\eeqn
%%
%where $\bar N(\Theta,\delta)$ is the tight entropy  analogue of the Sobolev covering numbers.  These lower bounds for ${\bf AC}$ are our vehicle for proving expectation lower bounds.
 In most examples, there is a critical $\eta=\eta(n,\Theta)$ after which \eref{accf} decreases exponentially.  This critical value $\eta(\Theta, n)$ is essential since it  yields, as a consequence  bounds  of type $ e_m(\Theta,\hf)\le C\eta(\Theta, n)^q$ which have been extensively studied in statistics, for
 \beqn
\label{exp}
e_m(\Theta,\hf)=\sup_{\rho\in\cM(\Theta)}\E_{\rhom}\|f_z-f_\rho\|^q
\eeqn
 
 To evaluate lower bounds for the function ${\bf AC}_m(\Theta,f_\bz,\eta)$, \cite{dkpt} considered :
 
$${\bf AC}_n(\Theta,\eta):=\inf _{{ \hat f}}\sup_{{\rho}\in\cM(\Theta)}{\rhon}\{\|{ f_\rho}-{ \hat f}\|>\eta\}$$

and the following result has been established:

\begin{eqnarray*}
\label{s1}
{\bf AC}_n(\Theta,\eta)\ge C\{
\begin{array}{ll}
e^{-cn\eta^2}, \quad & \eta\ge \eta_n,\\
 1,& \eta\le \eta_n,  
\end{array}
\end{eqnarray*}
where $\eta_n$ is defined by the relation :
 $\ln {\bar N(\Theta,\eta_n)}\sim c^2n(\eta_n)^2$. $N(\Theta,\eta_n)$ is the 'tight entropy'  defined  by:
\begin{eqnarray*}
\bar N(\Theta,\eta): = sup \{N: \exists \ f_0, f_1,...f_N\in \Theta,~~\hbox{with}\; 
    c_0\eta
\leq   \|f_i-f_j\|_{L_2({ \rho_X})} \leq  c_1 \eta,\ \forall i\neq j\}.
\end{eqnarray*}
 % \begin{eqnarray*}
%\label{s1}
%\inf _{{ \hat f}}\sup_{{ \rho}\in\cM(\Theta)}{\rhon}\{\|{ f_\rho}-{ \hat f}\|>\eta\}\ge C\{
%\begin{array}{ll}
%e^{-cn\eta^2}, \quad & \eta\ge \eta_n,\\
 %1,& \eta\le \eta_n,  
%\end{array}
   %\end{eqnarray*}
For instance, 
 $\eta_n=n^{-\frac{s}{2s+d}}$  for the Besov space $B_q^s(L_\infty(\Rd))$
 which corresponds to similar results proved   in statistics (actually with more restricted assumptions on the set of probabilities $\rho$):
   $$\inf _{{ \hat f}}\sup_{{ \rho}\in\cM'(B_q^s(L_\infty(\Rd)))}\E\|{ f_\rho}-{ \hat f}\|_{dx}\ge c n^{-\frac{s}{2s+d}}.$$
 See, for instance \cite{IH},  \cite{Stone},  \cite{Nem} for a slightly more restricted context than Besov spaces, and \cite{djkp93}...

Concerning upper bounds for ${\bf AC}_n(\Theta,\eta)$, many reverse properties have been established:see for instance \cite{BY} in statistical context, \cite{CS},
\cite{dkpt},   \cite{KT}, in learning theory.
These upper bounds are generally proved using particular estimation methods more often based on empirical mean square minimization.
$${ \hat f}=Argmin\{\sum_{i=1}^n(Y_i-f(X_i))^2,\; f\in \cH_n\}$$
   These very nice estimation rules raise nevertheless two important problems :
 First, they generally require  heavy computation times. The second serious problem lies in the fact that   their construction (the choice of $\cH_n$) is most of the time highly depending on $\Theta$ :  There also exist universal estimates (see \cite{tem}), however these rules are up to now prohibitive in terms of computation time.
 
 Our aim in this paper will be to emphasize special constructions and cases where it is possible to provide very simple (in terms of computation) estimates  enjoying  in addition the property of being universal : their construction does not depend on a particular $\Theta$ and still they have almost optimal properties for a whole bunch of spaces $\Theta$.
 These estimates are constructed using a thresholding schema, which has proven in the last decade in statistics to have very good properties for recovering  signals with inhomogeneous smoothness.
 
 In this paper, we will basically consider two particular situations. In the first case, we consider the RKHS situation. In this case, we produce a new  algorithm and investigate   its performances in $L_2(\hat\rho_X)$. The exponential rates of convergences are good:
 the critical value $\eta_n$ is the one predicted by \cite{dkpt}, and the exponential rates are  comparable to those recently obtained by \cite{SZ05}, although the loss is not the same ($L_2(\rho_X)$ in SZ), and the regularity assumptions are somewhat different : in SZ, regularity assumptions are expressed in terms of RKHS spaces. These assumptions may seem more intrinsic. However it is difficult to figure out exactly what they mean since they are depending on the unknown measure $\rho_X$. Our conditions are also depending on the kernel, but easy to figure out. 
 
 The second case considers a more specified situation where the $X_i$'s are one dimensional and the estimator is a wavelet thresholding estimate. The results are comparable in this setting to those obtained in the RKHS situation as concern the critical value and the exponential rates. The advantage here is that we are able to state the results in the $L_2(\rho_X)$ norm and the regularity conditions are expressed in terms of standard H\"older spaces.

\section{Least squares and thresholding procedures}
 
In this short section, we will consider the construction of our thresholding estimates. To make easier their understanding and motivate their consideration, we give here a connection with general least square estimates. However this construction   will not be used in the sequel and can be skipped  by a hurried reader which can go directly to the next section.

Empirical mean square minimization consists in considering
$${ \hat f}=Argmin\{\sum_{i=1}^n(Y_i-f(X_i))^2,\; f\in \cH_n\}$$
for a specified set $\cH_n$.
Let us look at particular cases of $\cH_n$ leading to
especially computable forms of $\hf$. Let us suppose that we have a collection of functions $(e_k)_k$ verifying the following property :
$${(P) : (e_k) :\frac1n\sum_{i=1}^ne_k(X_i)e_l(X_i)=\delta_{kl} }$$
 (i.e. $(e_k)$ is an orthonormal system for the empirical measure $\hat\rho$ on the $X_i's$.)
 $$\hat \rho_X=\frac1n\sum_{i=1}^n\delta(X_i)$$
 if $\delta(x)$ is the Dirac measure at the point $x$.
 
Now, associated to this collection of functions, let us consider the following particular spaces :
$$\cH_n^{(1)}=\{f=\sum_{i=1}^N\alpha_ie_i\},\quad
 \cH_n^{(2)}=\{f=\sum_{i=1}^N\alpha_ie_i,\; \sum|\alpha_i|\le \kappa\}
 $$
 $$ \cH_n^{(3)}=\{f=\sum_{i=1}^N\alpha_ie_i,\;\#\{|\alpha_i|\not=0\}\le \kappa\} 
$$
If we now introduce the 3 following  estimations of these coefficients :
$$\hat\alpha_k = \frac1n\sum_{i=1}^n e_k(X_i)Y_i, \quad  \hat\alpha_k^{(1)}= sign ({\hat\alpha_k})| \hat\alpha_k-\lambda|_+$$
$$  \hat\alpha_k^{(2)}=\hat\alpha_k I\{| \hat\alpha_k|\ge \lambda\}
$$
It is easy to prove that  there exists $\lam^i(\kappa)$ such that the following rules are empirical minimizers for the respective spaces $\cH_n^{(i)},\; i\in\,\{1,2,3\}$:
$$ { \hat f}^1=\sum_{k=1}^N {\hat\alpha_k}e_k\label{f1},\quad
{ \hat f}^2=\sum_{k=1}^N{ \hat\alpha_k^{(1)}}e_k$$
$$
{ \hat f}^3=\sum_{k=1}^N { \hat\alpha_k^{(2)}}e_k\label{f3}
$$

These three rules are common in the statistical litterature. $\hf^1$ is generally refered to as linear estimate, whereas, $\hf^2$ and $\hf^3$ are known as (respectively soft and hard) thresholding estimates.

Our aim in this paper is to study the behavior of these estimators, principally $\hf^3$, in  different situations.
The main difficulty of this paradigm obviously lies in the question :
 {How to choose the functions $(e_k)$  such that  condition $(P) $ is verified and suitably chosen tuning constants $ N,\; \lambda$ ?}
 
This first problem is difficult to solve, if not impossible, and in the sequel, we will not assume that property $(P)$ is verified,  but we are going to consider situations where this property can be considered as 'almost true'.

\section{RKHS situation}
\subsection{Assumptions, estimation rules and regularity conditions}
\subsubsection{Assumptions on the kernel}

Let us take the case of a symmetric  kernel  $K(\cdot,\cdot)$ (we do not explicitely need the fact that $K$ is a Mercer kernel).
We  assume that the kernel $K$ is uniformly bounded by an absolute constant $\kappa$.
Our fundamental assumption will be the following : 

(A) :
There exists a set of $p$ determinist points in $\R^d$ 
$$\{x_1,\ldots x_p\}$$
 ($p$ will tend to infinity with $n$) such that the following $p\times p$ matrix $M_{np}$ whose entries are, 
$(\frac1n \sum_{i=1}^nK(x_l,X_i)K(X_i,x_k))_{kl}$ is almost diagonal, in the sense that :
There exists $0\le\delta<1$ such that : 
\begin{equation}\label{cond}
\forall x\in \R^p,\;     \       \|x\|_{{ l}_2}^2(1-\delta)^2\le x^t\Mnp x \leq   \|x\|_{{ l}_2}^2(1+\delta)^2  
\end{equation}
\begin{equation}\label{cond1}
\|x\|_{{ l}_\infty}(1-\delta)\le\norm{\Mnp x}_{{ l}_\infty}
\end{equation}

%On the other side , it is clear that 
%\begin{equation}
 %u^t\Mnp u  \leq  \|u\|_{{ l}_2}^2\tau^2
%\end{equation}

We do not assume anything about $\delta$ but this quantity will enter into the performances results of the procedure. $\delta$ will be desired to be as small as possible.
Notice that in general, such an assumption reflects the concentration properties of the kernel, and is quite easy to verify in practical situations where $\delta$ can be computed empirically. In particular, we allow in the sequel $\delta$ to be a random quantity depending on the observations.

\subsubsection{Estimation rule}
 
Let us consider the following estimation rule :
We will denote by $Y$ the vector with coordinates $Y_i$,  $\eps_i=Y_i-f_\rho(X_i)$, and $\eps$ will be the vector with coordinates $\eps_i$.
Let us denote by $f_X$ the $n$ dimensional vector which entries are $f_\rho(X_i)$, and $K$ the $p\times n$ matrix which entries $ K(x_l,X_i)$ (so $ \frac 1n  K K^t = M_{np}$), and  introduce :
\begin{eqnarray}
 t_n &=&\frac{\log n}n,\quad
\lan=T \sqrt{t_n },
\\ z&=&(z_1,\ldots,z_p)^t=( K K^t )^{-1} K  Y,
\\ \tilde z &=&(\tilde z_1,\ldots,\tilde z_p)^t,\quad
\tilde z_l = z_lI\{|z_l|\ge \lan\}
\end{eqnarray}
%JE CROIS 
%$$\lan=T t_n ,\;  t_n=\sqrt{\frac{\log n}n}$$

$T$ will be chosen so that $T>\sqrt{M^2+\frac12}\vee 4$,
and finally, our estimate will be :
\begin{equation}\label{estimerc}
\hf=\sum_{l=1}^p\tilde z_lK(x_l,\cdot).
\end{equation}
As is easily seen, $\hat f$ takes its inspiration into $\hat f_3$ and it is worthwhile to notice that its construction do not depend on any regularity parameter.

\subsubsection{Regularity conditions}
We will assume the following sparsity conditions on the function $f_\rho$:

Let us take $$p=\PE{\left(\frac n{ \log n}\right)^{\frac {1}{2}} }$$

For any $n$, there exists $\alpha_1,\ldots,\alpha_p$, such that
\begin{eqnarray}
\|f_\rho-\sum_{l=1}^p\alpha_lK(x_l,\dot)\|_\infty&\le& cp^{-3/2}\label{Binfty}   \\
%n^{-1/2}\\
\forall \lam>0,\; card\{|\alpha_l|\ge \lam\}&\le&  c \lam^{-\frac 2{1+2s}} \label{weak}
\end{eqnarray}

These conditions reflect approximation properties for the function $f_\rho$ by linear combinations of vectors in the RKHS (when $K$ is a Mercer kernel).
These properties are quantified by conditions on the coefficients $\alpha_i$'s, which are standard in various situations (Fourier,  wavelet coefficients...).
As discussed in \cite{kepi00} condition (\ref{Binfty}) reflects a 'minimal compacity condition' which do not interfere in the entropy calculations (for instance) neither in the minimax rates of convergence. Condition  (\ref{weak}) does drive the rates. It is given here 
with a Lorentz type constraint on the $\alpha_i$'s. These conditions are obviously implied by $l_r$ conditions (for appropriate $r$) which then looks very much like Besov conditions.

We will measure the error by the following norm (empirical norm) :
\begin{equation}\label{norm}
\|g\|_{\hat \rho_X}^2=\frac1n\sum_{i=1}^ng(X_i)^2
\end{equation}
 $\PE{x}$ denotes the integer part of $x$.
 Our result is the following :
\begin{theorem}\label{mercer} 
Let us take $$p=\PE{\left(\frac n{ \log n}\right)^{\frac {1}{2}} }$$
For any $s>1/2$, we define,
 $$\eta_n=[\frac n{ \log n}]^{\frac {-s}{1+2s}}.$$
Under the conditions above,  there exists a constant $D$, such that
\begin{equation}
\sup_{{\rho}\in\cM(\Theta)}{\rhon}\{\|{ f_\rho}-{ \hat f}\|_{\hat \rho}>(1-\delta)^{-1}\eta\}\le
 T\{
\begin{array}{ll}
e^{-\gamma[n{ p^{-1}}\eta^2 \vee \log n]}, \; & \eta\ge D\eta_n,\\
 1,& \eta\le D\eta_n,  
\end{array}
\end{equation}

%as long as
% $$ [\frac n{\log n}]^{\frac {1}{1+2s}}\le p\le [\frac n{{ \log n}}]^{\frac {1}{2}}$$

\end{theorem}
\begin{rem}
As mentioned in the introduction these results prove that the behavior of this estimator is optimal in terms of  the critical value $\eta_n$ as predicted in \cite{dkpt}. In terms of exponential rates, they are suboptimal because of the term $p^{-1}$. However it is worthwhile to notice that these rates still are good: 
they are  comparable to those obtained by \cite{SZ05}, although the loss is not the same and the regularity assumptions are somewhat different.
In addition, we observe that if not entirely opimal, these rates are always better than $n^{-c}$.\\
 Finally, it is important to notice the following technical facts which will be crucial in the sequel : because $s>1/2$,
$\eta_n  \geq p^{-1}.$
Condition (\ref{weak}) can obviously always be replaced by:
\begin{equation}
\forall \lam>0,\; card\{|\alpha_l|\ge \lam\}  \le   c\lam^{-\frac 2{1+2s}}\wedge p) 
\end{equation}

\end{rem}
\subsection{Proof of the theorem}

First, let us remark that :

\begin{eqnarray*}
\|{ f_\rho}-{ \hat f}\|_{\hat \rho}&\le&\|f-\sum_{l=1}^p\alpha_lK(x_l,\dot)\|_\infty
+\|\sum_{l=1}^p\alpha_lK(x_l,\dot)-{ \hat f}\|_{\hat \rho}\\ &\le&
cp^{-3/2}+\|\sum_{l=1}^p(\alpha_l-\tilde z_l )K(x_l,\dot)\|_{\hat \rho}
\\ &\le&
cp^{-3/2}+ [(\alpha-\tilde z )M_{np} (\alpha-\tilde z )]^{\frac12}
\\ &\le& c\eta_n+ (1+\delta)[\sum_{l=1}^p(\alpha_l-\tilde z_l )^2]^{\frac12}
 \end{eqnarray*}
 Notice that the first line used hypothesis (\ref{Binfty}), and the last one (\ref{cond}).

\begin{eqnarray*}
\sum_{l=1}^p(\alpha_l-\tilde z_l )^2
&\le& \sum_{l=1}^p(\alpha_l- z_l )^2  \I\{|z_l|\ge \lan\}[\I\{|\alpha_l|\ge \lan/2\}+\I\{|\alpha_l|< \lan/2\}]\\
&+&
\sum_{l=1}^p[\alpha_l]^2\I\{|z_l|< \lan\}[\I\{|\alpha_l|\ge 2\lan\}+\I\{|\alpha_l|< 2\lan\}]
\\
&:=& BB+BS+SB+SS
\end{eqnarray*}
Let us study the term $SS$. First we remark that because of condition (\ref{weak}) on $\f$, we know that 
$$
  card\{|\alpha_l|\ge\lan\}\le c\lan^{\frac{-2}{1+2s}}\label{p*}
$$
and it is not difficult to prove that \eref{weak} is equivalent to the following characterization (the result is standard in Lorenz spaces and in any case can be found in \cite{cdkp}
\begin{equation}
\label{weak*}
\forall\; \lam>0, \sum_l\alpha_l^2\I\{|\alpha_l|< \lam\}\le c\lam^{\frac{4s}{1+2s}}
\end{equation}

%then, if $\ka\ge M/2$,  $\bjk|\ge \lan/2$ implies $j\le j_s$, 
Hence, using \eref{weak*} :
\\
\begin{eqnarray*}
SS&\le& c\lan^{\frac{4s}{1+2s}}=c(\sqrt{ t_n} T)^{\frac{4s}{1+2s}}  =cT^{\frac{4s}{1+2s}}  \eta_n^2
\end{eqnarray*}
Let us now investigate the term SB : We observe that $\I\{|z_l|< \lan\}\I\{|\alpha_l|\ge 2\lan\}\le \I\{|\alpha_l- z_l|\ge |\alpha_l|/2\}\I\{|\alpha_l|\ge 2\lan\}$, hence :
\begin{eqnarray*}
SB&\le& \sum_{l=1}^p[\alpha_l]^2 \I\{|z_l-\alpha_l|\ge |\alpha_l|/2\}\I\{|\alpha_l|\ge 2\lan\}\\
&\le&4\sum_{l=1}^p(\alpha_l- z_l )^2\I\{|\alpha_l|\ge 2\lan\}\\
 \end{eqnarray*}
 In the same way :
 \begin{eqnarray*}
BB&   =   & \sum_{l=1}^p[\alpha_l -z_l]^2 \I\{|z_l|\ge \lan  ;   \    |\alpha_l|\ge \lan/2\}
\\
&\le& \sum_{l=1}^p(\alpha_l- z_l )^2\I\{|\alpha_l|\ge \lan/2 \}
\\
 \end{eqnarray*}
  So $BB$ and $SB$ can be treated in the same way,
 since
  $$\sum_{l=1}^p(\alpha_l- z_l )^2\I\{|\alpha_l|\ge 2\lan\} \leq 
  \sum_{l=1}^p(\alpha_l- z_l )^2\I\{|\alpha_l|\ge \lan/2 \} $$
  $$BB + SB \leq  5 \sum_{l=1}^p(\alpha_l- z_l )^2\I\{|\alpha_l|\ge \lan/2 \} $$
  Let
  \begin{equation}
p^*= card\{|\alpha_l|\ge\lan/2 \}\le c(\lan/2)^{\frac{-2}{1+2s}}\label{p*}
\end{equation}
 
\subsubsection{ Study of   \protect$\sum_{l=1}^p(\alpha_l- z_l )^2\I\{|\alpha_l|\ge \lan/2 \}  $}

Let us denote by
$\bar f_X$  the vector with coordinates $[\bar f_X]_i=\bar f(X_i) =\sum_{l=1}^p\alpha_lK(x_l,X_i)$:
$$\bar f_X =   K^t \alpha$$
Let us recall that $f_X$ is  the $n$ dimensional vector which entries are $f(X_i)$.
and by hypothesis (\ref{Binfty}),  $|f(X_i)-\bar f(X_i) | \leq c p^{-3/2}$
So that, 
\begin{eqnarray*}
\alpha &=& (K K^t)^{-1} K\bar f_X,\\
z & =&   (K K^t)^{-1}K  Y=   (K K^t)^{-1}K [  f_X+\eps],\\
\alpha -z &=&  (K K^t)^{-1} K \epsilon     +    (K K^t)^{-1} K [\bar f_X-f_X] 
\end{eqnarray*}
>From this we deduce,
 $$\|\alpha-z\|_{l_2} \le\|(K K^t)^{-1} K \epsilon  \|_{l_2(p)}+ \|(K K^t)^{-1} K [\bar f_X-f_X]\|_{l_2(p)}$$
But, since $(K K^t)^{-1} =  \frac 1n \Mnp^{-1}$, and using (\ref{cond}),
\begin{eqnarray}\nonumber
 \|(K K^t)^{-1} K [\bar f_X-f_X]\|_{l_2(p)} &=&\frac 1n  \|M_{np}^{-1}K [\bar f_X-f_X]\|_{l_2(p)} \\ \nonumber&\leq&
 (1-\delta)^{-1}    \|  \frac 1n K [\bar f_X-f_X] \|_{l_2(p)} \\\nonumber
&\leq& (1-\delta)^{-1}     \kappa\|f_X-\bar f_X\|_{\infty}
\sqrt p \\ &\leq&  (1-\delta)^{-1} c \frac 1p \kappa
 \leq c (1-\delta)^{-1}   \kappa  \eta_n\label{21}
\end{eqnarray}
%$$ \|(K K^t)^{-1} K [\bar f_X-f_X]\|_{l_2(p)}^2 = $$
%$$ \|\Mnp^{-1}K_x[\bar f_X-f_X]\|_{l_2} \leq   \|(K K^t)^{-1} K[\bar f_X-f_X]\|_{l_2}   $$
%and$$\| K_x[u]\|^2_{l_2} =\sum_{l=1}^p (\frac 1{\sqrt n} \sum_{i=1}^n K(x_l,X_i)u_i)^2 $$
 %$$\leq     \|u\|_{l_2(n)}^2 \frac 1n \sum_{l=1}^p \sum_{i=1}^n K^2(x_l,X_i)=\|u\|_{l_2(n)}^2trace(M_{np})$$
% $$\leq p (1+\delta)\|u\|_{l_2(n)}^2 $$
%( ou bien $$\leq \tau^2 p (\frac 1{\sqrt n} \sum_{i=1}^n  |u_i|)^2\leq \tau^2  p    \|u\|_2^2 $$
%de toute facon ici PB !)

% So  \begin{eqnarray}
%\|\alpha-z\|_{l_2}&\le&\|\Mnp^{-1}K_x\eps\|_{l_2}+ \|\Mnp^{-1}K_x[\bar f_X-f_X]\|_{l_2}\nonumber \\
%&\le&\|\Mnp^{-1}K_x\eps\|_{l_2}+(1+\delta)(1-\delta)^{-1} p^{1/2} \sqrt n    p^{-5/2}
%\nonumber\\
%&\le&\|\Mnp^{-1}K_x\eps\|_{l_2}+(1+\delta)(1-\delta)^{-1}\eta_n\label{alphas}
%\end{eqnarray}
>From the calculations above and \eref{cond}, we deduce,

\begin{eqnarray}\label{pass}
\sum_{l=1}^p(\alpha_l&-&z_l)^2\I\{|\alpha_l|\ge \lan/2\}
\le
\sum_{l=1}^p ( (K K^t)^{-1} K \epsilon )_l^2        \I\{|\alpha_l|\ge \lan/2\}+ c(1-\delta)^{-1}[   \kappa\eta_n]^2
%&\le&(1-\delta)^{-2}\sum_l\I\{|\alpha_l|\ge \lan/2\}[\frac1n \sum_i\eps_iK(x_l,X_i)]^2
%+[c(1-\delta)^{-1}\eta_n]^2
\end{eqnarray}
Let us now recall the following inequality due to Pinelis \cite{Pin}, assuming that the  $\xi_i$'s are Hilbert space valued, independent random variables, such that $ \|\xi_i  -\E(\xi_i) \| \le \tilde M$ and 
$\E \|\xi_i  -\E(\xi_i) \|^2\le \si^2(\xi)$, 
\beqn 
Prob\big(\norm{\frac1n\sum_{i=1}^n[\xi_i-\E\xi_i]}\ge \lam\big)\le
2\exp\big\{\frac{-n\lam^2}{ 2( \lambda \tilde M/3+ \si^2(\xi) )  }\big\}\label{pine}
\eeqn
Now as $  \si^2(\xi) \leq \tilde M^2$, replacing $\si^2(\xi)$ in the RHS, we get:
\beqn 
Prob\big(\norm{\frac1n\sum_{i=1}^n[\xi_i-\E\xi_i]}\ge \lam\big)\le
2\exp\big\{\frac{-n\lam^2}{ 2( \lambda \tilde M/3+ \tilde{M}^2)  }\big\}
     \label{pine1}
\eeqn
As only $\lambda \leq \tilde M  $ is significant, since $Prob\big(\norm{\frac1n\sum_{i=1}^n[\xi_i-\E\xi_i]}\ge \lam\big)=0$, for $\lam>\tilde M$, 
\beqn 
Prob\big(\norm{\frac1n\sum_{i=1}^n[\xi_i-\E\xi_i]}\ge \lam\big)\le
2\exp\ \frac{-3n\lam^2}{ 8 \tilde{M}^2  } 
    \label{pine2}
\eeqn

Let us now take $\xi_i \in \R^p :$
$$(\xi_i)_l=(K(x_l,X_i)\eps_i)_{l}$$
in such a way that,
 $$\sum_{i} \xi_i = K\epsilon$$
  and the $\xi_i$ are independent.
It is easy to verify that $\E (\xi_i) =0.$
\\
Let us for all $U\in \R^p$ define the following Hilbertian norm :
$$\|U\|_A^2 = \sum_{l=1}^p ( n ( K K^t)^{-1}U)_l^2\I \ \{|\alpha_l|\ge \lan/2\} 
=\sum_{l=1}^p (  ( M_{np})^{-1}U)_l^2\I \ \{|\alpha_l|\ge \lan/2\}  $$
Then, 
$$\sum_{l=1}^p ( (K K^t)^{-1} K \epsilon )_l^2       \I \ \{|\alpha_l|\ge \lan/2\} = \|\frac 1n \sum_{i} \xi_i \|_A^2$$
Now, we have using (\ref{cond1})
\begin{eqnarray*}
 \|  \xi_i  \|_A^2 &=&\sum_{l=1}^p (   M_{np}^{-1}\xi_i)_l^2\I \ \{|\alpha_l|\ge \lan/2\}
 \\
 & \leq& p^*  (\sup_l (   M_{np}^{-1}\xi_i)_l)^2  \\
 &\leq& p^*  (\sup_l (  K(x_l,X_i)\epsilon_i)^2 \frac 1{(1-\delta)^2}
 \\
& \leq& p^*  (\kappa \epsilon_i)^2 \frac 1{(1-\delta)^2} \leq p^*  \frac{(M \kappa)^2}{(1-\delta)^2} 
\end{eqnarray*}

Now,  using \eref{pine2},   
$$Prob\big(\norm{\frac1n\sum_{i=1}^n[\xi_i-\E\xi_i]}^2 \ge \frac{(a \eta)^2}{(1-\delta)^2}\big)
%\le  2\exp\big\{   \frac{-n\frac{(a \eta)^2}{(1-\delta)^2}   }
%{ \frac{(a \eta)}{(1-\delta)}\frac 23 \sqrt{ p^*}  \frac{(M \kappa)}{(1-\delta)}+  2p^*  \frac{(M \kappa)^2}{(1-\delta)^2) } }\big\}  $$
%$$ \leq 2\exp\big\{   \frac{-n a^2 \eta^2   }{ \frac 23 a \eta  \sqrt{ p^*}   M \kappa   + 2 p^* M^2 \kappa^2 }\big\}
\leq 2\exp - \{ \frac 38 n\eta^2 \frac{a^2}{p^*M^2 \kappa^2}\} $$
So for $a>0$ suitably chosen, and taking account that $\eta >\eta_n$
%Now, we are using \eref{pine}, 
%$\xi_i=(K(x_l,X_i)\eps_i)_{l}$, noting that $\|\xi_i\|\le M\kappa (p^*)^{\frac12}$,
%$\E \|\xi_i\|\le (p^*)M^2\kappa^2$, and $\E\sum\xi^2\le (p^*)M\kappa/n$, for
%$\eta^2\ge (p^*) M\kappa/n$ (which is the case since hence since $(p^*)/n\le  ct_n^{\frac{4s}{1+2s}}<<\eta_n^2$.)

%$$\rhon(\sum_l\I\{|\alpha_l|\ge \lan/2\}[\frac1n \sum_i\eps_iK(x_l,X_i)]^2
%\ge \frac{(c \eta)^2}{(1-\delta)^2)\le \exp\frac{-n\eta^2}{2(p^*)M^2\kappa^2+M\kappa (p^*)^{\frac12}\eta}$$
%Hence, if $(p^*)^{\frac12}\eta\le (p^*)$, i.e. $\eta\le (p^*)^{\frac12}$, we have :
%$$\rhon(\sum_l[\frac1n \sum_i\eps_iK(x_l,X_i)]^2\ge \eta^2/2)\le \exp -c[n\eta^2p^{-1}\vee \log n]$$
%Because of the condition $s>1/2$.

\begin{eqnarray*}
\rhon&(&\sum_{l=1}^p(\alpha_l-z_l)^2\I\{|\alpha_l|\ge \lan/2\} \geq \frac{(2a \eta)^2}{(1-\delta)^2}\big)\\
&\le&
\rhon( \sum_{l=1}^p ( (K K^t)^{-1} K \epsilon )_l^2        \I\{|\alpha_l|\ge \lan/2\}+[ c \kappa(1-\delta)^{-1}\eta_n]^2\geq \frac{(2a \eta)^2}{(1-\delta)^2}\big)
\\ & \leq &\rhon(\sum_{l=1}^p ( (K K^t)^{-1} K \epsilon )_l^2        \I\{|\alpha_l|\ge \lan/2\}  \ge \frac{(a \eta)^2}{(1-\delta)^2}\big)
%\leq 2\exp\big\{   \frac{-n a^2 \eta^2   }{ \frac 23 a \eta  \sqrt{ p^*}   M \kappa   + 2 p^* M^2 \kappa^2 }\big\} $$
 \leq 2\exp - \{ \frac 38 n\eta^2 \frac{a^2}{p^*M^2 \kappa^2}\} 
\end{eqnarray*}

Now, if we recall that
$\eta_n=( \frac{\log n}{n})^{s/(1+2s)} ;       \     p^{-1} = \sqrt{t_n} ;    \   \lan=T\sqrt{t_n} $ and $p^* \le 4c (Tt_n)^{\frac{-1}{1+2s}} \wedge p       $, evaluation at the point
% if $(p^*)^{\frac12}\eta\le (p^*)$, i.e. $\eta\le (p^*)^{\frac12}$, 
 $\eta =\eta_n$ gives :
$$ 2\exp - \{ \frac 38 n\eta_n^2 \frac{c^2}{p^*M^2 \kappa^2}\} =  2\exp - \{ \frac 38 \log n \frac{a^2}
{c(\frac 1T)^{2/1+2s}M^2 \kappa^2}\}.$$
Hence
$$\rhon(\sum_l[\frac1n \sum_i\eps_iK(x_l,X_i)]^2\ge \eta^2/2)\le \exp -C[n\eta^2p^{-1}\vee \log n  ]$$

 \subsubsection{ Study of $ \sum_{l=1}^p(\alpha_l- z_l )^2\I\{|z_l|\ge \lan\}\I\{|\alpha_l|< \lan/2\}$    }

It remains now to study the term:
$$BS= \sum_{l=1}^p(\alpha_l- z_l )^2\I\{|z_l|\ge \lan\}\I\{|\alpha_l|< \lan/2\}\le\sum_{l=1}^p(\alpha_l- z_l )^2\I\{|z_l-\alpha_l|\ge \lan/2\}$$
Using the previous result with $p$ instead of $p^*$, we get,
\begin{eqnarray*}
\rhon &(&BS \geq \frac{(2a \eta)^2}{(1-\delta)^2})\\
& \leq  & \rhon(\sum_{l=1}^p(\alpha_l-z_l)^2  \geq \frac{(2a \eta)^2}{(1-\delta)^2}\big)\\
 &\leq &
 \rhon(\sum_{l=1}^p ( (K K^t)^{-1} K \epsilon )_l^2       \ge \frac{(a \eta)^2}{(1-\delta)^2}\big)
%\leq 2\exp\big\{   \frac{-n a^2 \eta^2   }{ \frac 23 a \eta  \sqrt{ p^*}   M \kappa   + 2 p^* M^2 \kappa^2 }\big\} $$
 \\
 &\leq &2\exp - \{ \frac 38 n\eta^2 \frac{a^2}{pM^2 \kappa^2}\} 
\end{eqnarray*}

We proceed as in the previous section, and obtain using (\ref{cond1}) :
\begin{eqnarray*}
|\alpha_l-z_l|&\le& |[(KK^t)^{-1}K\eps]_l|+ \|(KK^t)^{-1}K[\bar f_X-f_X]\|_{l_\infty}\nonumber \\
&\le& |\frac 1n (\Mnp^{-1}K\eps)_l|+(1-\delta)^{-1}  \kappa   p^{-3/2}
\nonumber\\
&\le&(1-\delta)^{-1}\|\frac 1n K\eps\|_{l_\infty}+(1-\delta)^{-1}\sqrt{t_n} M\kappa
\end{eqnarray*}

So 

$$\rhon(\exists l\in\;\{1,\ldots,p\},\; |\alpha_l-z_l|\ge \lan)
\leq \rhon(  (1-\delta)^{-1}  \sup_{l}|\frac 1n (K\eps)_l|+(1-\delta)^{-1} \kappa  p^{-3/2}\ge \lan)$$
$$ \leq  \sum_{l=1}^p \rhon(  (1-\delta)^{-1}| \frac 1n\sum_{i}K(x_l,X_i)\epsilon_i |+(1-\delta)^{-1}\kappa  p^{-3/2}\ge  T \sqrt{t_n})$$
$$ \leq  \sum_{l=1}^p \rhon(   | \frac 1n\sum_{i}K(x_l,X_i)\epsilon_i | \ge \sqrt{ \log n/n} (\frac{T}{1-\delta}
-\kappa ( \frac{\log n}{n})^{1/4} )$$

But for $n$ large enough 
$$(\frac{T}{1-\delta}
-\kappa ( \frac{\log n}{n})^{1/4} ) \geq  \frac{T}{2(1-\delta)}$$
and using Hoeffding inequality 
$$ \rhon(   | \frac 1n\sum_{i}K(x_l,X_i)\epsilon_i |  \geq  \frac{T}{2(1-\delta) \sqrt{ \log n/n }}) 
\leq  2 \exp - \frac{  T^2 \log n }{ 8 (1-\delta)^2 \kappa^2M^2}  $$
So
$$\rhon(\exists l\in\;\{1,\ldots,p\},\; |\alpha_l-z_l|\ge \lan) \leq 2p  \exp - \frac{  T^2 \log n }{ 8 (1-\delta)^2 \kappa^2M^2} \leq  Cn^{-\alpha}$$
with $\alpha >0$ if $T$ is large enough.

So :
$$\rhon( BS \geq  \frac{ a \eta}{1-\delta})  \leq 2\exp - \{ \frac 38 n\eta^2 \frac{a^2}{pM^2 \kappa^2}\}\wedge n^{-\alpha} $$

 This   yields the results.

\section{Wavelet results}
\subsection{Assumptions and estimation rules}
\subsubsection{Assumptions on the model}

In this section, we will concentrate on the case of dimension 1: the random variables $X_i$'s are now taking their values in $\X=$ compact domain of $\R$.  This case can easily be generalized to the case where the measure $\rho_X$ is a tensor product of measures
$\rho_{X_i}$, $i=1,\ldots, d$. However the full generalization to dimension $d$ is more involved and will not be discussed in this paper. In the case $d=1$, we  define the distribution function
$G$  such that $$ \forall  t  \in \R,   ~   G(t)=\rho(X\le t)   \in [0,1]$$
 and assume that it is a derivable function. We also define, 
$$ \forall x \in [0,1] ,   ~  { G^{-1}(x)=\inf\{t\in \R,\;  G(t)\ge x\}} .$$

Again, we will assume that $\f$ has sparsity conditions which can be in this case directly expressed in terms of regularity conditions. More precisely, we will denote by
$\M(\Theta_s)$, the set of measures $\rho$ verifying all the assumptions above with in addition the fact that $\f(G^{-1})\in B_\infty^s(L_\infty([0,1]))(M)$ (the ball of radius $M$ of the Besov  space). Notice that as we will only consider the case where $s>0$ (in fact $s>1/2$) $\f$ will always be bounded by $M$.

% Let us recall the following facts which are true under classical properties of regularity and moment
%vanishing (see \cite{Meyer} ):
 
 Let us consider $\{\psi_{j,k}, j\ge \underline{j}+1, \; 0 \leq k <2^j\}$ a   wavelet
basis  on $[0,1]$ (at least continuously differentiable, with enough moment conditions; the  length of
 of the support of $\psi_{j,k}$ the will be supposed to be less than $N 2^{-j}$). We recall that : $\psi_{\underline{j},k}=\phi_{\underline{j}k}$ denotes the scaling function.
These assumptions are standard (see \cite{cdv}).

Let us expand $f$ in the wavelet basis : $$f (G^{-1})= \sum_{j= \underline{j}}^\infty  \sum_{k \in }\beta_{j,k} \psi_{j,k}.$$ and
 it  is well known that for $0 \leq \gamma <\infty,$
$ f (G^{-1})$ belongs to
$B^\gamma_{\infty}(L_\infty([0,1]))$ iff  (and we will take this as the  $B^\gamma_{\infty}(L_\infty([0,1]))-$norm): 
$$ \sup_{j\ge \underline{j}} 2^{j(\gamma +\frac12 )} \sup_{ 0\leq k <2^j }|\beta_{j,k}|  =:
\|f\|_{B^\gamma_{\infty}}
 < \infty .$$
 %This is precisely the norm that will be taken here.

In this section  our loss will be  measured in term of $\L_2
$, with respect to the measure $d\rho_X$:
$$ \|f\|_{{ \rho_X}}=[\int f(x)^2d\rho_X(x)]^{\frac 12}.$$

\subsubsection{Estimation Algorithm :}
Again, we put 
$$ t_n:=\frac{\log n}n,\quad
\lan=\kappa \sqrt{t_n },
$$
define :
$$\hat G_n(x)=\frac1n\sum_{i=1}^nI\{X_i\le x\},$$
and let us  introduce the ordered statistic : $X_{(1)}\le\ldots\le X_{(n)}$. Doing this, we introduce a new ordering on the indices $\{1,\ldots,n\}$. Keeping this ordering, we denote $Y_{(1)},\ldots, Y_{(n)}$. Note that $Y_{(1)},\ldots Y_{(n)}$ is generally not the ordered statistic of $Y_{1},\ldots Y_{n}$.
\\
The estimator is constructed in the following way:

\begin{itemize}

 \item Step 1: Estimation of the wavelet coefficients : $${\hbjk}=\frac1n\sumin Y_i\psijk(\hat G_n(\frac in))=\frac1n\sumin Y_{(i)}\psijk(X_{(i)})$$
  \item 
Step 2: Thresholding $$ {\tbjk }= \hbjk \I\{|\hbjk|\ge \lambda_n\}$$
 \item Step 3: Reconstruction
$${ \hat f}=\sum_{j= \underline{j}}^{{J}}\sum_k
{\tbjk}\psijk(\hat G_n)$$
\end{itemize}

Note that this  algorithm is an adaptation of the standard wavelet algorithm introduced in \cite{dj-a} in the case of an equispaced design. It has been investigated in \cite{KP04}, where the expectation properties of the $L_p(dx)$ losses have investigated (instead of here the deviation properties of the $\L_2 d\rho_X
$). It proves to have very powerful properties.
One of them is its remarkable simplicity in terms of computation.
To illustrate this, we give here the main steps of the computation algorithm :
\\
{\it Algorithm :
\begin{enumerate}   \item  
Sort the $X_i$'s, \item Change the numbering in such a way that $X_i$ has rank
$i$, \item Calculate the highest level $alpha$-coefficients using the formula :
$$
\hat \alpha_{J'k}=\frac{1}{n}\sum_{i=1}^n\phi_{J'k}( i/n)Y_i,\quad
(2^{J'}=n)$$
\item Calculate the wavelet coefficients using the classical pyramidal algorithm
\item Perform a thresholding algorithm giving rise to $\tilde
\beta_{jk}$ coefficients,\item  Reconstruct the estimator, using again the
standard backward pyramidal algorithm,  obtaining
$$\hat f=\sum_{j= \underline{j}}^J\sum_{ 0\leq k <2^j }\tilde \beta_{jk}\psi_{jk}(\hat
G_n(x))$$
which is a function especially easy to draw.
\end{enumerate}}

Our aim in this section is to prove the following theorem.
\begin{theorem}
With the conditions above,
 $\forall  s>\frac12$,
  $$\eta_n=[\frac n{ \log n}]^{\frac {-s}{1+2s}},$$
there exist positive constants $\gamma,\; T,\; D$ such that,
 \begin{eqnarray*}
\sup_{ \rho\in\M(\Theta_s)}{\rhon}\{\|{ f_\rho}-{ \hat f}\|>\eta\}\le T\{
\begin{array}{ll}
e^{-\gamma[n{ 2^{-J}J^{-1}}\eta^2 \vee \log n]}, \; & \eta\ge D\eta_n,\\
 1,& \eta\le D\eta_n,  
\end{array}
\end{eqnarray*}  
%\begin{eqnarray*}
%\sup_{ \rho\in\M(\Theta_s)}{\rhon}\{\|{ f_\rho}-{P_{V_J} \hat f}\|>\eta\}\le C\{
%\begin{array}{ll}
%e^{-cn{ 2^{-J}}\eta^2\vee\log n }, \; & \eta\ge \eta_n,\\
% 1,& \eta\le \eta_n,  
%\end{array}
%\end{eqnarray*}
as long as
 $$ [\frac n{\log n}]^{\frac {1}{1+2s}}\le 2^{ J}\le [\frac n{{ \log n}}]^{\frac {1}{2}}$$
 \end{theorem} 
 \begin{rem}
As mentioned in the introduction these  results are comparable  to those obtained in the RKHS situation as concern the critical value and the exponential rates. The advantage here is that we are able to state the results in the $L_2(\rho_X)$ norm and the regularity conditions are expressed in terms of standard H\"older spaces.
We expressed the results in a slightly different way, leaving the choice of $J$, as an option. If we optimize oour results in $J$, we take $ 2^{ J}= [\frac n{\log n}]^{\frac {1}{1+2s}} $ which gives better rate results but fails in being adaptive. If we want our estimate to be universal (work for any $s>1/2$) we need to take $2^{ J}\le [\frac n{{ \log n}}]^{\frac {1}{2}}$.

\end{rem}
 
\subsection{Proof of the theorem}

Throughout the proof, the constant $c$ will denote a constant which may vary from one line to the other, but may be explicitely calculated. For a sake of simplicity we will not make explicite the constants obtained in the proof (although it could be done easily) since we do not think that they are optimal in any sense.

It will be essential in the sequel to notice that with the assumptions above, we have :

$$ \|f\|_{{L_2( \X\rho_X)}}= \|f({  G^{-1}})\|_{L_2([0,1],dx)}.$$

 Since $ \|f\|_{{ \rho_X}}= \|f({  G^{-1}})\|_{dx}$, we have if
 $$f_\rho(G^{-1})=\sumjk\bjk\psijk$$
 \begin{eqnarray*}
\|\hf-f_\rho\|_{ \rho_X}&=& \|\hf(G^{-1})-f_\rho(  G^{-1})\|_{dx}
\\
&=&
\|\sum_{j= \underline{j}}^{J}\sum_k
{\tbjk}\psijk(\hat G_n(G^{-1}))-\sumjk\bjk\psijk\|_{dx}
\\
&\le & \|\sum_{j= \underline{j}}^{J}\sum_k
{\tbjk}[\psijk(\hat G_n(G^{-1}))-\psijk]\|_{dx}+ \|\sum_{j= \underline{j}}^{J}\sum_k
[{\tbjk}-\bjk]\psijk\|_{dx}\\ &+&\|\sum_{j=\ge J+1}\sum_k\bjk\psijk\|_{dx}
\end{eqnarray*}
Hence

\begin{eqnarray*}
\|\hf-\f\|_{ \rho_X}^2
&\le &3 
[\|\sum_{j= \underline{j}}^{J}\sum_k
{\tbjk}[\psijk(\hat G_n(G^{-1}))-\psijk]\|_{dx}^2+\sum_{j= \underline{j}}^{J}\sum_k[{\tbjk}-\bjk]^2+\sum_{j\ge J+1}\sum_k\bjk^2]
\\
&\le &
(I)+(II)+(III)
\end{eqnarray*}

If $\f(G^{-1})\in B_\infty^s(L_\infty([0,1]))(M)$, then 
 $$III =\sum_{j\ge J+1}\sum_k\bjk^2 \leq  \sum_{j\ge J+1} 2^j \sup_{k}\bjk^2 \leq 
M^2  \sum_{j\ge J+1} 2^j 2^{-j(2s+1)} \leq  M^2 2^{-2Js}\leq M^2 \eta_n^2$$
 if $2^J\ge t_n^{\frac {-1}{1+2s}} = (\sqrt{\frac{\log n}n})^{\frac {-1}{1+2s}} $
 
    Let us now study the second term:
\begin{eqnarray*}
(II)&\le& \sum_{j= \underline{j}}^{J}\sum_k[{\hbjk}-\bjk]^2\I\{|\hbjk|\ge \lan\} \left[   \I\{|\bjk|\ge \lan/2\}+\I\{|\bjk|< \lan/2\}  \right]  \\
&+&
\sum_{j= \underline{j}}^{J}\sum_k[\bjk]^2\I\{|\hbjk|< \lan\}    \left[ \I\{|\bjk|\ge 2\lan\}+\I\{|\bjk|< 2\lan\}
       \right]\
\\
&:=& BB+BS+SB+SS
\end{eqnarray*}
Let us study the term $SS$. First we remark that, as $\f(G^{-1})\in B_\infty^s(L_\infty([0,1]))(M)$, then $|\bjk|\le M2^{-j(s+\frac12)}$, hence if we denote :
$$2^{j_s}=t_n^{\frac{-1}{1+2s}}$$
%then, if $\ka\ge M/2$,  $\bjk|\ge \lan/2$ implies $j\le j_s$, hence :
 \begin{eqnarray*}
SS&\le& \sum_{j= \underline{j}}^{j_s}\sum_k[\bjk]^2\I\{|\bjk|<2\lan\}+\sum_{j=j_s}^{J}\sum_k[\bjk]^2\I\{|\bjk|<2\lan\}\\
&\le&\sum_{j= \underline{j}}^{j_s}\sum_k[2\lan]^2+\sum_{j=j_s}^{J}\sum_k[\bjk]^2\\
\\
&\le& 22^{j_s}(2\lan)^2+\sum_{j=j_s}^{J}2^jM^22^{-2j(s+\frac12)}\\
&\le&(8\kappa^2+2M^2)\eta_n^2 
\end{eqnarray*}
Let us now investigate the term SB : We observe that $$\I\{|\hbjk|< \lan\}\I\{|\bjk|\ge 2\lan\}\le \I\{|\hbjk-\bjk|\ge |\bjk|/2\}\I\{|\bjk|\ge 2\lan\}$$, hence :
\begin{eqnarray*}
SB&\le& \sum_{j= \underline{j}}^{J}\sum_k[\bjk]^2 \I\{|\hbjk-\bjk|\ge |\bjk|/2\}\I\{|\bjk|\ge 2\lan\}\\
&\le&4\sum_{j= \underline{j}}^{J}\sum_k|\hbjk-\bjk|^2\I\{|\bjk|\ge 2\lan\}\\
%&\le&4BB
\end{eqnarray*}
       So
       $$BB +SB \leq 5\sum_{j= \underline{j}}^{J}\sum_k|\hbjk-\bjk|^2\I\{|\bjk|\ge \lan/2\} =5 BB'$$
Now, we investigate the term $BB'$.

If we recall that $X_{(1)}\le\ldots\le X_{(n)}$. Doing this we introduce a new ordering on the indices $\{1,\ldots,n\}$, and that we keep this ordering, to denote $Y_{(1)},\ldots, Y_{(n)}$. We also introduce $U_i=G(X_i),\; i=1,\ldots,n$,
as well as the associated $U_{(1)},\ldots, U_{(n)}$. Notice  that the $U_{(i)}$'s are ordered (since $G$ is increasing) and the $U_i$'s are i.i.d. uniformly distributed. 
\begin{eqnarray}
\hbjk-\bjk&= &\frac1n\sumin Y_{(i)}\psijk(\frac in)-\bjk\nonumber\\
&= &[\frac1n\sumin \f(G^{-1}(U_{(i)}))\psijk(\frac in)-\int\psijk \f(G^{-1})]+
[\frac1n\sumin\eps_{(i)}\psijk(\frac in)]\nonumber
\\
&= &
[\frac1n\sumin \f(G^{-1}(U_{(i)}))\psijk(\frac in)-\int\psijk \f(G^{-1})]+
[\frac1n\sumin\eps_{i}\psijk(\frac in)]\nonumber
\\
&:= & A_{jk}+B_{jk}\label{num}
\end{eqnarray}
%We have used here the independence  between the $\eps_i$'s and the $X_i$'s to replace $\eps_{(i)}$ by $\eps_i$ (the equality is true in distribution). Notice that this could be untrue, if this assumption was removed.  
%Let us, for convenience  put $U_{(0)}=0$ and $U_{(n+1)}=1$
Let us begin by the following lemma which proof is obvious (but which will be  useful in the sequel :
\begin{lemma}\label{lem}
For any $r\ge 1$, we have
\begin{equation}
\frac 1n\sum_{i=1}^n|\psijk(\frac in)|^r\le \tau_r2^{j(\frac r2-1)}+\tau_r'\frac{2^{j(1+\frac r2)}}n \label{lem}
\end{equation}
with $\tau_r=N\|\psi\|_\infty$ and $\tau_r'=Nr\|\psi'\|_\infty(\|\psi\|_\infty)^{r-1}$
\end{lemma}

\noindent Let us put:
$$\hat F_n(x)=\frac1n\sumin\I\{U_{i}\le x\},\quad\Delta_n:=\sup_{x\in [0,1]}|\hat F_n(x)-x|. $$
and $\bar s=s\wedge 1$, using (\ref{lem}) for the third inequality, 
\begin{eqnarray}\nonumber
|A_{jk}|&\le&\frac1n\sumin |\f(G^{-1}(U_{(i)}))-\f(G^{-1}(\frac in))||\psijk(\frac in)|
\\ &+&\sum_{i=1}^n\int_{(i-1)/n}^{i/n}|\f(G^{-1}(x)\psijk(x)-\f(G^{-1}(\frac in))\psijk(\frac in)|\nonumber
\\
&\le&\Delta_n^{\bar s}\|\f(G^{-1})\|_{\bar s \infty\infty}\frac1n\sumin|\psijk(\frac in)|\nonumber
\\&+&\sum_{i=1}^n\int_{(i-1)/n}^{i/n}[2^{j/2}\|\psi\|_\infty\|\f(G^{-1})\|_{\bar s \infty\infty}n^{-\bar s}        +\|\psi\|_{1 \infty\infty}\|\f(G^{-1})\|_\infty\frac{2^{3j/2}}n]\I\{x\in[\frac k{2^j},\frac{ k+N}{2^j}]\}dx
 \nonumber
\\
%&\le&\Delta_n^{\bar s}\|f(G^{-1})\|_{\bar %s\infty\infty}[\int|\psijk|+\frac1n\sumin|\psijk(\frac in)|-\int|\psijk|]
%\\&+&[ \|\psi\|_\infty\|f(G^{-1})\|_{\bar s \infty\infty}n^{-\bar s}2^{\frac %j2}+\|\psi\|_{1 \infty\infty}\|f(G^{-1})\|_\infty\frac{2^{3j/2}}n]\frac1n\sumin[\I\{\frac %in\in[\frac k{2^j},\frac{ k+N}{2^j}]\}+\I\{\frac{ i-1}n\in[\frac k{2^j},\frac{ %k+N}{2^j}]\}]\\
&\le&\Delta_n^{\bar s}\|\f(G^{-1})\|_{\bar s \infty\infty}\{\tau_12^{-j/2}+\tau_1'\frac{2^{3j/2}}n\}\nonumber
\\ &+& N\|\psi\|_\infty\|\f(G^{-1})\|_{\bar s \infty\infty}n^{-\bar s}2^{-\frac j2}
+N\|\psi\|_{1  \infty\infty}\|\f(G^{-1})\|_\infty\frac{2^{j/2}}n\nonumber
\\
&\le& C_1\Delta_n^{\bar s}2^{-j/2}+C_2\frac{2^{j/2}}n+C_3n^{-\bar s}2^{-\frac j2}\label{num*}
\end{eqnarray}
where
$$C_1=\tau_1+\tau_1',\; C_2=N\|\psi'\|_\infty,\; C_3=N\|\psi\|_\infty \|\f(G^{-1})\|_{\bar s \infty\infty}$$

The last line uses the fact that for $j\le J , ~  2^{2j} \le n$.
We can then state the following lemma :
\begin{lemma}
For $J$ such that $t_n^{\frac{-1}{1+2s}}\le 2^J\le t_n^{-1/2}$, we have :
\begin{equation}\label{Ajk}
\rhon(\sum_{j= \underline{j}}^{J}\sum_kA_{jk}^2\ge \eta^2)\le \exp-Cn2^{-J}\eta^2 \vee\log n,
\end{equation}
for all $\eta\ge D\eta_n$, where $C=2(2C_1^2N)^{\frac{-1}{\bar s}}$
\end{lemma}

\noindent
{\sl Proof of the lemma :}

We observe that $$\sum_{j= \underline{j}}^{J}\sum_k[\frac{2^{j/2}}n]^2\le c\frac{2^{2J}}{n^2}\le c\frac1n<<\eta_n^2 .$$
  $$ \hbox{and} ~~ \sum_{j= \underline{j}}^{J}\sum_kn^{-2\bar s}2^{-j}\le Jn^{-2\bar s}<<\eta_n^2$$. 
  \begin{equation}
\sum_{j= \underline{j}}^{J}\sum_k\Delta_n^{2\bar s}2^{-j}\le   J\Delta_n^{2\bar s}\label{num2}
\end{equation}
 Hence,
 for $\eta\ge D\eta_n$, and $n$ large enough,
\begin{eqnarray*}
\rhon(\sum_{j= \underline{j}}^{J}\sum_kA_{jk}^2\ge \eta^2)&\le&\rhon(C_1^2J\Delta_n^{2\bar s}\ge \eta^2/2)\\
&\le&K\exp{-cn[\eta J^{-1/2}]^{\frac2{\bar s}}}\I\{\eta^2  \le 2C_1^2J\}
\end{eqnarray*}
The last line uses the following Dvoreski, Kiefer and Wolfovitz bound 
 (see for instance the review on the
subject in Devroye Lugosi section 12.) : For any $ \lambda>0$, there
exists a  universal constant
$K$, such that:

\begin{eqnarray}  
\P(   \Delta_n\ge \lambda)\le K\exp-2n\lambda^2 
 \label{DKW1}
\end{eqnarray} 
(and noticing that $\Delta_n\le 1$)
\\
Now, for $s\ge 1$, $n[\eta]^{\frac2{\bar s}}J^{\frac{-1}{2\bar s}}=n\eta^2J^{-1/2}\ge
n\eta^22^{-J} \vee\log n $. \\
Identically, for $1/2<s<1$ and $\eta\ge D\eta_n$,
\begin{eqnarray*}
n[\eta]^{\frac2{ s}}J^{\frac{-1}{2 s}}&\ge& n\eta^22^{-J}\eta_n^{2(\frac1s-1)}2^JJ^{\frac{-1}{2 s}}\\ 
&\ge&n\eta^22^{-J}2^{-2sj_s(\frac1s-1)}2^{j_s}J^{\frac{-1}{2 s}}
\\ 
&\ge&n\eta^22^{-J}2^{j_s(2s-1)}J^{\frac{-1}{2 s}}\ge n\eta^22^{-J} \vee\log n
\end{eqnarray*}
This ends up the proof of the lemma.\hskip 6cm \qed
\\
Let us now investigate the term  corresponding to the $B_{jk}$'s. We have the following lemma :

\begin{lemma}\label{lem1}
For $J$ such that $t_n^{\frac{-1}{1+2s}}\le 2^{J}\le t_n^{-1/2}$, there exists a constant $c$ such that:
\begin{equation}\label{Bjk}
\rhon(\sum_{j= \underline{j}}^{J}\sum_kB_{jk}^2\I\{|\bjk|\ge \lan/2\}\ge \eta^2)\le c\exp-cn2^{-J}\eta^2 \vee\log n,
\end{equation}
for all $1\ge\eta\ge D\eta_n$
\end{lemma}

\noindent
{\sl Proof of the lemma :}

Let us first remark that since $\f(G^{-1})\in B_\infty^s(L_\infty([0,1]))(M)$, then $|\bjk|\le M2^{-j(s+\frac12)}$ and then, if $\ka\ge 2M$,  $|\bjk|\ge \lan/2$ implies $j\le j_s$, hence :
\begin{eqnarray*}
\sum_{j= \underline{j}}^{J}\sum_kB_{jk}^2\I\{|\bjk|\ge \lan/2\}&\le& \sum_{j= \underline{j}}^{j_s}\sum_kB_{jk}^2
\\&\le& \sum_{j= \underline{j}}^{j_s}2^j\sup_{k}B_{jk}^2\le 22^{j_s}\sup_{jk}B_{jk}^2
\end{eqnarray*}

We will investigate separately the cases $\eta\le 1$, and $\eta\ge1$.
Let us begin with the fist case:

\begin{eqnarray}
\rhon(\sum_{j= \underline{j}}^{J}
\sum_kB_{jk}^2
\I\{|\bjk|\ge \lan/2\}
\ge \eta^2)
&\le &
  \rhon(  2^{j_s+1}\sup_{jk} B_{jk}^2\ge \eta^2) \nonumber
\\
&\le& \sum_{j= \underline{j}}^{J}\sum_k \rhon(|\sumin\psijk(\frac in)\eps_i|\ge  n \eta 2^{-j_s/2}/\sqrt2)   \nonumber
\\
&\le&\sum_{j= \underline{j}}^{J}\sum_k\exp\{\frac{-n^2\eta^22^{-j_s}/2}{2(nC_3+n \eta M\|\psi\|_\infty 2^{(j-j_s)/2}/3) }\}\nonumber
 \\
&\le&2^{j_s+1}\exp\{\frac{-n\eta^2 2^{-j_s}}{4(C_3+ \eta M\|\psi\|_\infty /3) }\} \label{num3}
\end{eqnarray}
In the last line we used Bernstein inequality (cf Bernstein \cite{Bern}), since the variables
$\psijk(\frac in)\eps_i$ are a sequence of  independent  bounded random variables (by $M\|\psi\|_\infty 2^{\frac j2}$), with zero mean and 
 $$\E[ \sum_i \psijk(\frac in)\eps_i]^2 \leq C_3n$$
 ($M^2( \tau_2+\tau_2') := C_3$ using (\ref{lem}).)
\\
 Hence we obtain :
\begin{eqnarray}
\rhon(\sum_{j= \underline{j}}^{J}\sum_kB_{jk}^2\I\{|\bjk|\ge \lan/2\}\ge \eta^2)&\le&
2\exp\{-cn\eta^22^{-j_s}+ \frac 12\log n\}\label{bern}
\end{eqnarray}
with $c=4(C_3+M\|\psi\|_\infty/3)^{-1}$ since $\eta\le 1$. As $\eta\ge D\eta_n$, it is easy to see that for $D$ large enough,
$cn\eta^22^{-j_s}\ge 2\log n.$
Hence in this case, we get the bound : $\exp-cn2^{-J}\eta^2 \vee\log n$

%Hence, if $\nu_n$ is the empirical measure associated to the $\eps_i$'x : i.e. 
%$$
%\nu_n=\frac1n\sumin\delta_{\eps_i}
%$$
%what we have to investigate is the following probability :
%$$
%\rhon(\sup_{f\in \cF}\int fd\nu_n\ge c\eta)
%$$
%for
%$$\cF=\{f:\{1,\ldots,n\}\mapsto \R/\ f(i)=

Let us now study the case where $\eta\ge 1$, we'll use Mac Diarmid's inequality (see \cite{MacD} we have the following lemma :

\begin{lemma}
For $J$ such that $t_n^{\frac{-1}{1+2s}}\le 2^{J}\le t_n^{-1/2}$, we have :
\begin{equation}\label{Bjk1}
\rhon(\sum_{j= \underline{j}}^{J}\sum_kB_{jk}^2\I\{|\bjk|\ge \lan/2\}\ge \eta^2)\le \exp-Cn2^{-J}\eta^2 \vee\log n,
\end{equation}
for all $\eta\ge 1$, and $C= \frac 1{2B^2},\; B^2=2M^2N^2\|\psi\|_\infty$.
\end{lemma}

\noindent{\it Proof of the lemma :}
We have :
$$\rhon(\sum_{j= \underline{j}}^{J}\sum_kB_{jk}^2\I\{|\bjk|\ge \lan/2\}\ge \eta^2)\leq \rhon(F(\eps_1,\ldots,\eps_n)\ge \eta^2)$$
with:
 $$F(\eps_1,\ldots,\eps_l,\ldots,\eps_n) = \sum_{j= \underline{j}}^{j_s}\sum_k  \frac1{n^2}[\sumin\psijk(\frac in)\eps_i]^2
$$
\begin{eqnarray*}
|\Delta F_l|&=&|F(\eps_1,\ldots,\eps_l,\ldots,\eps_n)-F(\eps_1,\ldots,\eps_l',\ldots,\eps_n)|\\ &=&
\sum_{j= \underline{j}}^{j_s}\sum_k\I\{|\bjk|\ge \lan/2\}\frac1{n^2}\left ([\sumin\psijk(\frac in)\eps_i]^2-
[\sumin\psijk(\frac in)\eps_i+\psijk(\frac ln)(\eps_l'-\eps_l)]^2\right)
\\ &\le&2M^2\sum_{j= \underline{j}}^{j_s}\sum_{k,\; |\frac ln-\frac k{2^j}|\le \frac N{2^j}}\frac1{n^2}
\sumin|\psijk(\frac in)||\psijk(\frac ln)|
\\ &\le&2M^2N^2\frac1{n^2}\|\psi\|_\infty^2\sum_{j= \underline{j}}^{j_s}2^j\frac n{2^j}
\\ &\le&2M^2N^2\|\psi\|_\infty^2\frac Jn=:B^2 \frac{j_s}n
\end{eqnarray*}
On the other hand,
\begin{eqnarray*}
\E_\rhon F(\eps_1,\ldots,\eps_n) &\le&\sum_{j= \underline{j}}^{j_s}\sum_k\frac1{n^2}[\sumin\psijk(\frac in)\eps_i]^2
\\ &\le& \sum_{j= \underline{j}}^{j_s}\sum_k \frac1{n^2}\sumin\psijk(\frac in)^2M^2
\\ &\le& M^2C_3\sum_{j= \underline{j}}^{j_s}\sum_k \frac1{n}
\\ &\le& 2M^2C_3\frac{2^{j_s}}{n}\leq c\eta_n^2
\end{eqnarray*}
Hence, for $\eta\ge D\eta_n$,
\begin{eqnarray*}
\rhon(\sum_{j= \underline{j}}^{J}\sum_kB_{jk}^2\I\{|\bjk|\ge \lan/2\}\ge \eta^2)&\le&
\rhon(|F(\eps_1,\ldots,\eps_n)-\E_\rhon F(\eps_1,\ldots,\eps_n)| \ge \eta^2/2)
\\ &\le& \exp \frac {-2\eta^4}{4n(\frac {BJ}n)^2}\le \exp-nC\frac{\eta^4}{J^2}
\end{eqnarray*}

Now, for $\eta\ge 1$, we obviously have  $Cn\frac{\eta^4}{J^2}\ge Cn2^{-J}\eta^2\vee\log n$, which proves the result of the lemma.
\qed

Notice also that, using exactly the same proof, we have also the following result, which will be used later :
\begin{lemma}\label{bis}
For $J$ such that $t_n^{\frac{-1}{1+2s}}\le 2^J\le t_n^{-1/2}$, we have :
\begin{equation}\label{Bjk2}
\rhon(\sum_{j= \underline{j}}^{J}\sum_kB_{jk}^2\ge \lambda^2)\le \exp-Cn\lambda^4/J^2 \vee\log n,
\end{equation}
for all $\lambda^2\ge 2M^2C_1t_n^{1/2}$, $C= \frac 1{2B^2},\; B^2=2M^2N^2\|\psi\|_\infty$.
\end{lemma}

This achieves bounding the term (BB).
We now proceed to bound the term (BS):
  
\begin{eqnarray*}
\sum_{j= \underline{j}}^{J}\sum_k(\hbjk-\bjk)^2\I\{|\bjk|<\lan/2\}\I\{|\hbjk-\bjk|\ge\lan/2\}&\le&
\sum_{j= \underline{j}}^{J}\sum_k(\hbjk-\bjk)^2\I\{|\hbjk-\bjk|\ge\lan/2\}\\&\le&
2^{J+1}\sup_{jk}\{(\hbjk-\bjk)^2; \; |\hbjk-\bjk|\ge\lan/2\}
\end{eqnarray*}
Hence
\begin{eqnarray*}
\rhon(
\sum_{j= \underline{j}}^{J}\sum_k(\hbjk-\bjk)^2\I\{|\hbjk-\bjk|\ge\lan/2\}\ge \eta^2\})&\le&
2^{J+1}\rhon(|\hbjk-\bjk|\ge \eta 2^{-J/2}/2\vee \lan/2)
\end{eqnarray*}
Now, using (\ref{num}) and (\ref{num2}), we get
\begin{eqnarray*}
2^{J+1}\rhon(|\hbjk-\bjk|\ge \eta 2^{-J/2}/2&\vee&\lan/2)\le
2^{J+1}\rhon(C_1\Delta_n^{\bar s}2^{-J/2}\ge \eta 2^{-J/2}/4\vee \lan/8)\\&+&
2^{J+1}\rhon(|B_{jk}|\ge \eta 4^{-J/2}/2\vee \lan/8)\\ &\le&
2^{J+1}K\exp -cn(\frac\eta 4)^{\frac 2{\bar s}}\I\{\eta\le 4C_1\}
\\
&+&
2^{J+1}\exp\{\frac{-n(\eta^22^{-J}/16\vee \lan^2/64)}{2C_3+(\eta/4\vee \lan/8)2^{J/2}) M\|\psi\|_\infty}\}
\end{eqnarray*}
The first term may be bounded as in Lemma \ref{lem1}, the second one may be bounded by:\\
$ \exp-c[n2^J\eta^2 \vee\log  n],
$ with $c=(64C_3)^{-1}$
if $\eta\le 1$.
\\
For $\eta\ge 1$ , we have :
\begin{eqnarray*}
\sum_{j= \underline{j}}^{J}\sum_k(\hbjk-\bjk)^2\I\{|\bjk|<\lan/2\}&\le&
2^{J+1}\sup_{jk}A_{jk}^2+\sum_{j= \underline{j}}^{J}\sum_kB_{jk}^2
\end{eqnarray*}
Hence,
\begin{eqnarray*}
\rhon(
\sum_{j= \underline{j}}^{J}\sum_k(\hbjk-\bjk)^2\I\{|\bjk|<\lan/2\}\ge \eta^2)&\le&
\rhon(2^{J+1}\sup_{jk}A_{jk}^2\ge \eta^2/2)
\\ &+&\rhon(\sum_{j= \underline{j}}^{J}\sum_kB_{jk}^2\ge \eta^2/2\})
\end{eqnarray*}
The first term, treated as above, gives the same bound since in this case the condition $\eta\le 1$ was not necessary. For the second term, we use the lemma (\ref{bis}).

This achieves the proof for the term (II), which can be summarised in the following proposition.

\begin{proposition}

 $\forall  s>\frac12$
\begin{eqnarray*}
\sup_{ \rho\in\M(\Theta_s)}{\rhon}\{\sum_{j= \underline{j}}^{J}\sum_k[{\tbjk}-\bjk]^2>t^2\}\le c\{
\begin{array}{ll}
e^{-Cn{ 2^{-J}}t^2\vee\log n }, \; & t\ge \eta_n,\\
 1,& t\le \eta_n,  
\end{array}
\end{eqnarray*}
\\ if $ [\frac n{\log n}]^{\frac {1}{1+2s}}\le 2^{ J}\le [\frac n{{ \log n}}]^{\frac {1}{2}}$
for $C=(64C_3)^{-1}\wedge(2B^2)^{-1}\wedge 2(\frac2 {C_1N})^{\frac 1{\bar s}}$

\end{proposition}
It remains, now to study the term (I).

We have :

\begin{eqnarray*}
\|\sum_{j= \underline{j}}^{J}\sum_k
{\tbjk}[\psijk(\hat G_n(G^{-1}))-\psijk]\|_{dx}&\le&
\|\sum_{j= \underline{j}}^{J}\sum_k
{|\tbjk-\bjk|}[|\psijk(\hat G_n(G^{-1}))-\psijk]]\|_{dx}
\\
&+&\|\sum_{j= \underline{j}}^{J}\sum_k
{|\bjk|}[|\psijk(\hat G_n(G^{-1}))-\psijk]]\|_{dx}\\ &\le&
\|\sum_{j= \underline{j}}^{J}\sum_k
{|\hbjk-\bjk|}[|\psijk(\hat G_n(G^{-1}))-\psijk]]\|_{dx}
\\
&+&2\|\sum_{j= \underline{j}}^{J}\sum_k
{|\bjk|}[|\psijk(\hat G_n(G^{-1}))-\psijk]]\|_{dx}
\end{eqnarray*}
since $|\tbjk-\bjk|\le|\hbjk-\bjk|+|\bjk|$.
If $Z=\sum|\bjk|\psijk$ we observe that $\|Z\|_{s\infty\infty}=\|\f(G^{-1})\|_{s\infty\infty}$, so:
\begin{eqnarray*}
\|\sum_{j= \underline{j}}^{J}\sum_k
{|\bjk|}[|\psijk(\hat G_n(G^{-1}))-\psijk]]\|_{dx}
&=&\|Z(\hat G_n(G^{-1})-Z\|_{dx}
\\ &\le&
\|\f(G^{-1})\|_{\bar s\infty\infty}\Delta_n^{\bar s}
\end{eqnarray*}
Hence,

\begin{eqnarray*}
\rhon(\|\sum_{j= \underline{j}}^{J}\sum_k
{|\bjk|}[|\psijk(\hat G_n(G^{-1}))-\psijk]]\|_{dx}\ge\eta)&\le&
\rhon(\|\f(G^{-1})\|_{\bar s\infty\infty}\Delta_n^{\bar s}\ge\eta) \\ &\le&
K\exp{-cn\eta^{\frac 2{\bar s}}}\I\{\eta/\|\f(G^{-1})\|_{\bar s\infty\infty}\le 1\} \\ &\le&K  \exp\{-cn\eta^22^{-J}\vee\log n\}
\end{eqnarray*}

As above (see the proof of lemma \ref{lem1}), 
with $c=2\|\f(G^{-1})\|_{\bar s\infty\infty}^{\frac 2{\bar s}}$, here.
 
Concerning the stochastic term, using (\ref{num})  we have : 
\begin{eqnarray*}
\|\sum_{j= \underline{j}}^{J}\sum_k
{|\hbjk-\bjk|}[|\psijk(\hat G_n(G^{-1}))-\psijk]]\|_{dx}&\le&
\|\sum_{j= \underline{j}}^{J}\sum_k
{|A_{jk}|}[|\psijk(\hat G_n(G^{-1}))-\psijk]]\|_{dx}
\\ &+&
\|\sum_{j= \underline{j}}^{J}\sum_k
{| B_{jk}|}[|\psijk(\hat G_n(G^{-1}))-\psijk]]\|_{dx}
\end{eqnarray*}

Now, if
$Z'=\sum_{j= \underline{j}}^{J}\sum_k
{|A_{jk}|}\psijk$, using (\ref{num*}), and $s>\frac 12$,
\begin{eqnarray*}
\|Z'\|_{1/2\infty\infty}&\le&\sup_{j\le J,k}\{
2^j|A_{jk}|\}\\  &\le&\sup_{j\le J,k}\{C_1\Delta_n^{\bar s}2^{j/2}+C_2\frac{2^{3j/2}}n+C_3n^{-\bar s}\}\\ &\le& C_1\Delta_n^{\bar s}2^{J/2}+(C_2+C_3)\frac{2^{3J/2}}n\end{eqnarray*}

Let us investigate separately the two contributions:
As above,
\begin{eqnarray*}
\rhon(\|\sum_{j= \underline{j}}^{J}\sum_k
{|\hbjk-\bjk|}[|\psijk(\hat G_n(G^{-1}))-\psijk]]\|_{dx}\ge\eta)&\le&
\rhon(\|Z'\|_{1/2\infty\infty}\Delta_n^{1/2}\ge\eta) 
\end{eqnarray*}
Furthermore,
\begin{eqnarray*}
\rhon(\Delta_n^{1/2+\bar s}2^{J/2}\ge \eta/(2C_1))
\le\exp\{-2n(\frac\eta{2C_1}2^{-J/2})^{\frac2{\bar s+1/2}}\}\I\{\frac\eta{2C_1}2^{-J/2}\le 1\}
\end{eqnarray*}

Now, as $\bar s> 1/2$, we have, for $\eta2^{-J/2}\le 2C_1 $,
$n(\eta2^{-J/2})^{\frac2{\bar s+1/2}}\ge (2C_1)^{\frac{1-2\bar s}{1+\bar s/2}}n(\eta2^{-J/2})^2\vee\log n$, for $\eta\ge \eta_n$. 
 On the other hand, for $\tilde C= C_2+C_3$

\begin{eqnarray*}
\rhon(\tilde C\frac{2^{3J/2}}n \Delta_n^{1/2}\ge \eta)
\le\exp\{-n2(\tilde C)^{-4}(n\eta2^{-3J/2})^4\} \I\{n\eta2^{-3J/2}\le \tilde C\}
\end{eqnarray*}

And obviously, on the range we are considering $n(n\eta2^{-3J/2})^4\ge n(\eta2^{-J/2})^2\vee\log n$.

 Now for the last term, 
 ($\|\sum_{j= \underline{j}}^{J}\sum_k
{| B_{jk}|}[|\psijk(\hat G_n(G^{-1}))-\psijk]]\|_{dx}$), considering again the $U_{(i)}$'s and  putting $U_{(0)}=0,\; U_{(n+1)}=1$, we have, on $[U_{(i)}, U_{(i+1)}]$, $\hat G(G^{-1}(x))=\frac in$.
 For any arbitrary $a>0$, we have

\begin{eqnarray*}
\|\sum_{j= \underline{j}}^{J}\sum_k
{||}[|\psijk(\hat G_n(G^{-1}))-\psijk|]\|_{dx}^2
&\le&
[\sum_{j= \underline{j}}^{J}\|\sum_k
{|B_{jk}|}[|\psijk(\hat G_n(G^{-1}))-\psijk|]\|_{dx}]^2
\\&\le&
2^{(J+1)a}\sum_{j= \underline{j}}^{J}2^{-ja}\int[ \sum_k
{|B_{jk}|}|\psijk(\hat G_n(G^{-1}))-\psijk]]^2
\\ &\le&
2^{(J+1)a}\sum_{j= \underline{j}}^{J}2^{-ja}\sum_{i=0}^n\int_{U_{(i)}}^{U_{(i+1)}}[ \sum_k
{|B_{jk}|}|\psijk(\frac in)-\psijk(x)|]^2dx
\end{eqnarray*}
Now, we will distinguish two cases : either $\frac in \in [U_{(i)}-\frac N{2^j}, U_{(i+1)}+\frac N{2^j}]$ (case I) or not (case II, which implies that $\Delta_n2^j\ge N$).

In case I, if we denote by $\Delta_{n,i}=\sup\{|\frac in-U_{(i)}|,|\frac in-U_{(i+1)}|\}$, and $I_{jk}$ is the support of $\psijk$ ,as $\psi$ is continuously differentiable,
we get,  for $x\in [U_{(i)}, U_{(i+1)}]$ :

$[\sum_k
{|B_{jk}|}|
 \psijk(\frac in)-\psijk(x)]]^2\le [\sum_k
{|B_{jk}|}\|\psi'\|_\infty2^{3j/2}\Delta_{n,i}\I_{\Ijk}(x)]^2\le N\sum_k
{|B_{jk}|^2}\|\psi'\|_\infty^22^{3j}\Delta_{n,i}^2\I_{\Ijk}(x)$

The last inequality is true because only a finite number of $\I_{\Ijk}(x)$'s are not zero at the same time.

 If we now remark that in case I, $\Delta_{n,i}\le 2N2^{-j}\wedge \Delta_n$ we get,  for $x\in [U_{(i)}, U_{(i+1)}]$ :

$[\sum_k
{|B_{jk}|}|
 \psijk(\frac in)-\psijk]]^2(x) \le 2N^2\sum_k
{|B_{jk}|^2}\|\psi'\|_\infty2^{2j}\Delta_{n}\I_{\Ijk}(x).$
\\ \\
In case II, we get,  for $x\in [U_{(i)}, U_{(i+1)}]$, using again the fact that only a finite number of $\psijk$'s are not zero at the same time :

\begin{eqnarray*}
[\sum_k
{|B_{jk}|}|
 \psijk(\frac in)-\psijk(x)|]^2&\le& 2\left\{ [\sum_k|B_{jk}|
 |\psijk(\frac in)|]^2+ [\sum_k{|B_{jk}|}|
 \psijk(x)]]^2\right\}\I\{\Delta_n2^j\ge N \}
 \\
 &\le& 2\left[N\|\psi\|_\infty^22^j\sup_{j\le J,\; k} B_{jk}^2+[\sum_k
{|B_{jk}|}
| \psijk(x)|]^2\right]\I\{\Delta_n2^j\ge N \}
\end{eqnarray*}

 Putting the two cases together,  we deduce :
 
 \begin{eqnarray*}
\|\sum_{j= \underline{j}}^{J}&\sum_k&
{|B_{jk}|}[|\psijk(\hat G_n(G^{-1}))-\psijk|]\|_{dx}^2
\le
c2^{Ja}\sum_{j= \underline{j}}^{J}2^{-ja}\sum_{i=0}^n\int_{U_{(i)}}^{U_{(i+1)}}[ \sum_k
{| B_{jk}|}|\psijk(\frac in)-\psijk(x)|]^2dx
\\ &\le&
c2^{Ja}\sum_{j= \underline{j}}^{J}2^{-ja}\sum_{i=0}^n\int_{U_{(i)}}^{U_{(i+1)}}\left\{N^2\sum_k
{|B_{jk}|^2}\|\psi'\|_\infty2^{2j}\Delta_{n}\right.\\ &+&\left.\left[N\|\psi\|_\infty^22^j\sup_{j\le J,\; k} B_{jk}^2
+[\sum_k
{|B_{jk}|}
| \psijk(x)|]^2\right]\I\{\Delta_n2^j\ge N \}\right\}dx
 \\ &\le&
c2^{Ja}\sum_{j= \underline{j}}^{J}2^{-ja} \left\{N^2\sum_k
{|B_{jk}|^2}\|\psi'\|_\infty2^j\Delta_{n}\right.\\ &+&\left.\left[ N\|\psi\|_\infty^22^j\sup_{j\le J,\; k}B_{jk}^2+\sum_k
{|B_{jk}|}^2N^{-1}\Delta_{n}2^j\right]\I\{\Delta_n2^j\ge N \}\right\}
\\ &\le&
c\left[\sum_{j= \underline{j}}^{J}\sum_k
{|B_{jk}|^2}\Delta_{n}
 +2^j\sup_{j\le J,\; k} B_{jk}^2\I\{\Delta_n2^j\ge N \}\right]:= A+B
\end{eqnarray*}

To study the first term, again using lemma \ref{bis},  and (\ref{DKW1}), we get
\begin{eqnarray*}
\rhon(A\ge \eta^2/3)&\le&\rhon(\sum_{j= \underline{j}}^{J}\sum_{k}{| B_{jk}|}^22^{J}\Delta_n\ge c\eta^2)\\
&\le&\rhon(\sum_{j= \underline{j}}^{J}\sum_{k}{| B_{jk}|}^22^{J}\ge t^2)+\rhon( \Delta_n\ge c\eta^2/t^2)
\\&\le& c\exp -c[n\frac{t^4}{J^22^{2J}}\vee \log n]+K\exp -n \frac{c^22\eta^4}{t^4}
\end{eqnarray*}
for $t^22^{-J}\ge ct_n^{1/2}$ :
Optimizing in $t$, we find, for $t^4={c\eta^22^J J}$,
$$\rhon(A\ge \eta^2/3)\le \exp -n\eta^22^{-J}J^{-1}$$

This is valid if $t^22^{-J}\ge ct_n^{1/2}$ i.e. $\eta2^{-J/2}\ge cn^{-1/2}$.

Now taking $t=m J$, we find $$\rhon(A\ge \eta^2/3)\le\exp -[d\log n]$$
using again the fact that $s>\frac 12$ and $\eta\ge D\eta_n$.

On the other hand, we have also the following bound using Bernstein inequality (see (\ref{bern}):
\begin{equation}
\rhon(\sum_{j= \underline{j}}^{J}\sum_{k}{| B_{jk}|}^22^{j}\ge t^2)\le
2^J\rhon(| B_{jk}|^22^{2J}\ge t^2)\le
2^J\exp -n {ct^2}{2^{-2J}}
\end{equation}
For $t2^{-J/2}\le c'$.
If then again, we optimize in $t$, we find :
$t^2=\eta^{4/3}2^{2J/3}$ leading to the rate :
$\exp -n\eta^{4/3}2^{-4J/3}$
We have $\eta^{4/3}2^{-4J/3}\ge \eta^22^{-J}$ for $\eta\le 2^{-J/2}$.
In this case, we precisely have
$t^22^{-J/2}=\eta^{4/3}2^{2J/3}2^{-J/2}\le2^{-J/2}$.

It is obvious that the second  term ($B$ ) may be bounded (using (\ref{DKW1})) by
\begin{eqnarray*}
\rhon(\Delta_n2^j\ge N)\le K\exp -2nN^2 2^{-2J}\le \exp -2N\log n
\end{eqnarray*}

Now, we have, using (\ref{bern})
$$\rhon(\sup_{ \underline{j}\le j\le J,\; k} B_{jk}^22^J\ge c'\eta^2)\le c\exp cn\eta^22^{-J}$$
if $\eta\le c''$. Notice that the constant $c''$ may be chosen arbitrarily. Of course this choice will change the constant $c$. Hence, let us take $c''=MN$, and
now,
let us remark that,

$2^j B_{jk}^2\le 2^j[\frac 1n\sum_i M2^{J/2}\I\{\frac in\in [\frac k{2^j}, \frac {k+n}{2^j}]\}]^2\le 2^j[\frac 1n M2^{J/2}\frac {nN}{2^j}]^2\le M^2N^2$.
Hence the probability for $\sup_{ \underline{j}\le j\le J,\; k}B_{jk}^22^J$ to exceed $\eta^2$ is zero for $\eta^2> M^2N^2$.

This achieves bounding the term $SS$ as well as ends up the proof of the theorem.

\bibliographystyle{apalike}
\bibliography{Genbibl}

\def\cprime{$'$}
\begin{thebibliography}{}

\bibitem[Bernstein, 1946]{Bern}
Bernstein, S. (1946).
\newblock {\em The theory of Probability}.
\newblock Gastehizdal Publishing House, Moscow.

\bibitem[Cohen et~al., 1993]{cdv}
Cohen, A., Daubechies, I., and Vial, P. (1993).
\newblock Wavelets on the interval and fast wavelet transforms.
\newblock {\em Appl. Comput. Harmon. Anal.}, 1(1):54--81.

\bibitem[Cohen et~al., 2001]{cdkp}
Cohen, A., DeVore, R., Kerkyacharian, G., and Picard, D. (2001).
\newblock Maximal spaces with given rate of convergence for thresholding
  algorithms.
\newblock {\em Appl. Comput. Harmon. Anal.}, 11(2):167--191.

\bibitem[Cucker and Smale, 2002]{CS}
Cucker, F. and Smale, S. (2002).
\newblock On the mathematical foundations of learning.
\newblock {\em Bull. Amer. Math. Soc. (N.S.)}, 39(1):1--49 (electronic).

\bibitem[DeVore et~al., 2004]{dkpt}
DeVore, R., Kerkyacharian, G., Picard, D., and Temlyakov, V. (2004).
\newblock Mathematical methods for supervised learning.
\newblock Technical report, IMI.
\newblock University of South carolina.

\bibitem[Diarmid, 1989]{MacD}
Diarmid, M. (1989).
\newblock On the method of bounded differences.
\newblock In {\em Surveys in Combinatorics}, pages 148--188. Cambridge
  University Press, Cambridge.

\bibitem[Donoho and Johnstone, 1994]{dj-a}
Donoho, D.~L. and Johnstone, I.~M. (1994).
\newblock Ideal spatial adaptation by wavelet shrinkage.
\newblock {\em Biometrika}, 81(3):425--455.

\bibitem[Donoho et~al., 1995]{djkp93}
Donoho, D.~L., Johnstone, I.~M., Kerkyacharian, G., and Picard, D. (1995).
\newblock Wavelet shrinkage: Asymptopia?
\newblock {\em Journal of the Royal Statistical Society, Series B},
  57:301--369.
\newblock With Discussion.

\bibitem[Gy{\"o}rfi et~al., 2002]{4Authors}
Gy{\"o}rfi, L., Kohler, M., Krzy{\.z}ak, A., and Walk, H. (2002).
\newblock {\em A distribution-free theory of nonparametric regression}.
\newblock Springer Series in Statistics. Springer-Verlag, New York.

\bibitem[Ibragimov and Has'minski{\u\i}, 1981]{IH}
Ibragimov, I.~A. and Has'minski{\u\i}, R.~Z. (1981).
\newblock {\em Statistical estimation}.
\newblock Springer-Verlag, New York.
\newblock Asymptotic theory, Translated from the Russian by Samuel Kotz.

\bibitem[Kerkyacharian and Picard, 2000]{kepi00}
Kerkyacharian, G. and Picard, D. (2000).
\newblock Thresholding algorithms and well-concentrated bases.
\newblock {\em Test}, 9(2).

\bibitem[Kerkyacharian and Picard, 2004]{KP04}
Kerkyacharian, G. and Picard, D. (2004).
\newblock Regression in random design and warped wavelets.
\newblock {\em Bernoulli}, 10(6):1053--1105.

\bibitem[Konyagyn and Temlyakov, 2004]{KT}
Konyagyn, S.~V. and Temlyakov, V.~N. (2004).
\newblock Some error estimates in learning theory.
\newblock In {\em Approximation theory: a volume dedicated to Borislav
  Bojanov}, pages 126--144. Prof. M. Drinov Acad. Publ. House, Sofia.

\bibitem[Korostelev, 2003]{K}
Korostelev, A. (2003).
\newblock The {B}ahadur risk in probability density estimation.
\newblock {\em Statist. Decisions}, 21(2):139--148.

\bibitem[Korostelev and Spokoiny, 1996]{KS}
Korostelev, A.~P. and Spokoiny, V.~G. (1996).
\newblock Exact asymptotics of minimax {B}ahadur risk in {L}ipschitz
  regression.
\newblock {\em Statistics}, 28(1):13--24.

\bibitem[Nemirovskiy, 1985]{Nem}
Nemirovskiy, A.~S. (1985).
\newblock Nonparametric estimation of smooth regression functions.
\newblock {\em Izv. Akad. Nauk SSSR Tekhn. Kibernet.}, (3):50--60, 235.

\bibitem[Pinelis, 1994]{Pin}
Pinelis, I. (1994).
\newblock Optimum bounds for the distributions of martingales in banach spaces.
\newblock {\em Ann. Probab.}, 22:1679--1706.

\bibitem[Poggio and Smale, 2003]{PS}
Poggio, T. and Smale, S. (2003).
\newblock The mathematics of learning: dealing with data.
\newblock {\em Notices Amer. Math. Soc.}, 50(5):537--544.

\bibitem[Smale and Zhou, 2005]{SZ05}
Smale, S. and Zhou, D.-X. (2005).
\newblock Learning theory estimates via operators and their approximations.
\newblock Technical report, Toyota Technological Institute.

\bibitem[Stone, 1982]{Stone}
Stone, C.~J. (1982).
\newblock Optimal global rates of convergence for nonparametric regression.
\newblock {\em Ann. Statist.}, 10(4):1040--1053.

\bibitem[Temlyakov, 2005]{tem}
Temlyakov, V. (2005).
\newblock Approximation in learning theory.
\newblock Technical report, IMI.
\newblock University of South carolina.

\bibitem[Van~de Geer, 2001]{VG}
Van~de Geer, S. (2001).
\newblock {\em Empirical processes in M-estimation}.
\newblock Cambridge University Press, New York.

\bibitem[Yang and Barron, 1999]{BY}
Yang, Y. and Barron, A. (1999).
\newblock Information-theoretic determination of minimax rates of convergence.
\newblock {\em Ann. Statist.}, 27(5):1564--1599.

\end{thebibliography}

\end{document}